\documentclass[a4paper]{amsart}
\usepackage[utf8]{inputenc}
\usepackage[english]{babel}

\usepackage{amscd,amssymb}
\usepackage{amsmath}
\usepackage{amsthm}
\usepackage{graphicx}
\usepackage{fancybox}
\usepackage{epic,eepic}
\usepackage{amstext}
\usepackage{mathtools}
\usepackage{esint}
\usepackage[square,numbers]{natbib}
\usepackage{booktabs}
\usepackage[hide links]{hyperref}
\usepackage{comment}
\usepackage{mathtools}
\usepackage{tikz,pgfplots}
\usepackage{subcaption}

\usepackage[noabbrev, capitalise, nameinlink]{cleveref}
\crefname{equation}{}{}
\crefname{assumption}{Assumption}{Assumptions}
\crefformat{equation}{\textup{#2(#1)#3}}

\usepackage{algorithm}
\usepackage{algpseudocode}
\usepackage{tabularx}
\newcolumntype{Y}{>{\centering\arraybackslash}X}

\newtheorem{theorem}{Theorem}[section]
\newtheorem{corollary}[theorem]{Corollary}
\newtheorem{lemma}[theorem]{Lemma}

\theoremstyle{definition}

\theoremstyle{remark}
\newtheorem{remark}[theorem]{Remark}
\numberwithin{theorem}{section}
\numberwithin{equation}{section}
\numberwithin{figure}{section}

\def\Th{\mathcal{T}_h}

\def\pih{\pi_h}
\def\Pih{\Pi_h}
\def\Uh{U_h}
\def\Sh{\Sigma_h}

\def\diam{\operatorname{diam}}

\def\ddiv{\operatorname{div}}
\def\bu{B_u}

\def\jul{J_{u,\lambda}}
\def\with{\,:\,}
\def\dx{\,\text{d}x}
\DeclareMathOperator*{\argmin}{arg\,min}

\newcommand{\tn}[2]{{\left\lVert #1 \right\rVert}_{L^2(#2)}}
\newcommand{\tnf}[2]{\| #1 \|_{L^2(#2)}}

\newcommand{\tnof}[1]{\| #1 \|_{L^2}}

\newcommand{\tspof}[2]{{( #1,#2 )}_{L^2}}

\numberwithin{equation}{section}
\numberwithin{theorem}{section}

\AtBeginDocument{%
	\def\MR#1{}
}

\begin{document}
\title[Mixed finite elements for the Gross-Pitaevskii eigenvalue problem]{Mixed finite elements for the Gross-Pitaevskii eigenvalue problem: a priori error analysis and guaranteed lower energy bound}

\author[D.~Gallistl, M.~Hauck, Y.~Liang, D.~Peterseim]{Dietmar Gallistl$^*$, Moritz Hauck$^+$, Yizhou Liang$^\dagger$, Daniel Peterseim$^{\dagger,\ddagger}$}
\address{${}^*$ Institute of Mathematics, University of Jena,  Ernst-Abbe-Platz~2, 07743 Jena, Germany}
\email{dietmar.gallistl@uni-jena.de}
\address{${}^+$ Department of Mathematical Sciences, University of Gothenburg and Chalmers University of Technology, 41296 Göteborg, Sweden}
\email{hauck@chalmers.se}
\address{${}^{\dagger}$ Institute of Mathematics, University of Augsburg, Universit\"atsstr.~12a, 86159 Augsburg, Germany}
\address{${}^{\ddagger}$Centre for Advanced Analytics and Predictive Sciences (CAAPS), University of Augsburg, Universit\"atsstr.~12a, 86159 Augsburg, Germany}
\email{\{yizhou.liang, daniel.peterseim\}@uni-a.de}
\thanks{The work of D.~Gallistl, M.~Hauck, and D.~Peterseim is part of projects that have received funding from the European Research Council ERC under the European Union's Horizon 2020 research and innovation program
(DG: project DAFNE, grant agreement No.~891734; MH, DP: RandomMultiScales, grant agreement No.~865751). M.~Hauck is also supported by the Knut and Alice Wallenberg foundation postdoctoral program in mathematics
for researchers from outside Sweden (Grant No. KAW 2022.0260). Y.~Liang is supported by a Humboldt Research Fellowship for Postdocs from the Alexander von Humboldt Foundation. 
}

\begin{abstract}
We establish an a priori error analysis for the lowest-order Raviart-Thomas finite element discretisation of the nonlinear Gross-Pitaevskii eigenvalue problem. Optimal convergence rates are obtained for the primal and dual variables as well as for the eigenvalue and energy approximations. In contrast to conformal approaches, which naturally imply upper energy bounds, the proposed mixed discretisation provides a guaranteed and asymptotically exact lower bound for the ground state energy. The theoretical results are illustrated by a series of numerical experiments.
\end{abstract}



\maketitle

\section{Introduction}
We study the Gross-Pitaevskii equation (GPE), a nonlinear eigenvalue problem that describes the quantum states of bosonic particles at ultracold temperatures, so-called Bose-Einstein condensates. Given a convex Lipschitz domain $\Omega \subset \mathbb R^d$ ($d = 1,2,3$), the GPE searches for $L^2$-normalised eigenstates $\{u_j \with j \in \mathbb N\}\subset H^1_0(\Omega)$ and corresponding eigenvalues $\lambda_j \in \mathbb R$ such that
\begin{equation}
	\label{eq:gpe}
	-\Delta u_j + V u_j + \kappa|u_j|^2 u_j = \lambda_j u_j
\end{equation}
holds in the weak sense. Here $V\in L^\infty(\Omega)$ denotes a non-negative trapping potential that confines the particles to a particular region within the domain, and $\kappa$ is a positive constant. Note that all eigenvalues of~\cref{eq:gpe} are real and positive and that the smallest eigenvalue is simple. Assuming a non-decreasing ordering of the eigenvalues, this means that $0<\lambda_1< \lambda_2\leq \dots$.

The nonlinear eigenvalue problem~\cref{eq:gpe} is the Euler-Lagrange equation for critical points of the Gross-Pitaevskii energy
\begin{equation}
	\label{eq:energy}
	\mathcal{E}(v) \coloneqq \tfrac12\tspof{\nabla v}{\nabla v} + \tfrac12 \tspof{Vv}{v} + \tfrac{\kappa}{4} \tspof{|v|^2 v}{v},\quad v \in H^1_0(\Omega),
\end{equation}
subject to the $L^2$-normalisation constraint. Of particular physical interest is the ground state of the Gross-Pitaevskii energy, characterized by
\begin{equation}
	\label{eq:gs}
	u \in  \argmin_{v \in H^1_0(\Omega)\with \tnof{v} = 1} \mathcal E (v).
\end{equation}
We emphasize that under the above assumptions on $\Omega$ and $V$, the global energy minimiser exists and is unique up to sign. Furthermore, the ground state $u$ (up to sign) coincides with the eigenstate $u_1$  corresponding to the smallest eigenvalue $\lambda_1$ of~\cref{eq:gpe}. The minimal energy~$E$ is related to the smallest eigenvalue $\lambda_1$ by $\lambda_1 = 2E +\tfrac{\kappa}{2}\|u\|_{L^4}^4$.
Note that the above theoretical results on the Gross-Pitaevskii problem can be found, e.g., in~\cite{CCM10}.

There are a number of discretisations in the literature to approximate the ground state of the GPE. Such discretisations are typically based on $H^1_0$-conforming methods, such as continuous finite elements~\cite{Zho04,CCM10,CHZ11}, spectral and pseudospectral methods~\cite{CCM10,Bao13}, multiscale methods~\cite{HMP14b,MR4379968,henning2023optimal,PWZ23}, and mesh-adaptive methods~\cite{Dan10,He21}. Conforming discretisations have in common that the ground state energy is approximated from above, as the energy is minimised in a subspace. In this work we instead use a mixed finite element discretisation, which allows asymptotically exact lower bounds on the ground state energy. In the linear setting such an approach has recently been introduced in \cite{Gal23}.

In addition to the guaranteed lower energy bound, we provide a rigorous a priori error analysis of the proposed mixed finite element method for the GPE. We prove first-order convergence for the primal and dual variables in the $L^2$-norm and second-order convergence for the energy and eigenvalue approximations. So far, error estimates of this form have only been shown for conforming approximations of the primal variable. Although there is a large body of work on mixed discretisation methods for linear eigenvalue problems (see, e.g., the review article~\cite{Bof10}), nonlinear eigenvector problems have not yet been addressed. In fact, the present error analysis differs substantially from the established techniques used in the linear~case.

\section{Mixed finite element discretisation}
Consider a hierarchy of simplicial meshes $\{\Th\}_{h>0}$ of the domain $\Omega$, which we assume to be geometrically conforming (cf.~\cite[Def.~1.55]{ErG04})  and shape-regular (cf.~\cite[Def.~1.107]{ErG04}). We denote the elements of any mesh $\Th$ in the hierarchy by $K$ and define the mesh size  $h$ as the maximum diameter of elements in $\Th$, i.e., $h \coloneqq \max_{K \in \Th} \diam K$.
For the mixed discretisation of the Gross-Pitaevskii problem, we use the finite element pair $(\Sh,\Uh)$, where $\Sh$ denotes the lowest-order Raviart-Thomas finite element space with respect to~$\Th$ (see, e.g., \cite[Ch.~1.2.7]{ErG04}) and $\Uh$ is the space of $\Th$-piecewise constants. A discrete analogue of the gradient operator $G_h \colon \Uh \to \Sh$ is defined for arbitrary $v_h \in \Uh$ by the property 
\begin{equation}
	\label{eq:discgrad}
	\tspof{G_h v_h}{\tau_h} + \tspof{\ddiv \tau_h}{v_h} = 0
\end{equation}
for all $\tau_h\in\Sh$.
The discrete gradient gives rise to the discrete energy defined for any $v _h \in \Uh$ by 
\begin{equation*}
	\mathcal E_h(v_h) \coloneqq \tfrac12 \tspof{G_h v_h}{G_h v_h} + \tfrac12 \tspof{Vv_h}{v_h} + \tfrac{\kappa}{4} \tspof{|v_h|^2 v_h}{v_h}.
\end{equation*}
A discrete approximation $u_h \in \Uh$ of the ground state $u$ in the Raviart-Thomas space is then obtained as the solution of the finite-dimensional minimisation problem
\begin{equation}
	\label{eq:gsdisc}
	u_h \in \argmin_{v_h \in \Uh \with \tnof{v_h} = 1} \mathcal E_h(v_h).
\end{equation}
Note that in the discrete setting, the boundedness of the norms of the minimising sequence directly implies the strong convergence of a subsequence (Bolzano-Weierstrass theorem). Thus there always exist discrete energy minimisers $u_h$ and $-u_h$. Unlike in the continuous setting, cf.~\cref{eq:gs}, the solution $u_h$ to \cref{eq:gsdisc} is not unique up to sign in general. To have compatible signs of the ground state and its discrete approximation, we choose the sign of $u_h$ such that $(u,u_h)_{L^2}\geq 0$.

The proof of the guaranteed lower energy bound is based on certain properties of the operators $\pih\colon L^2(\Omega)\to \Uh$ and $\Pih\colon (L^2(\Omega))^d\to \Sigma_h$, which are defined as $L^2$-projections onto $\Uh$ and $\Sh$, respectively. By definition, $\pih$ and $\Pih$ are bounded with respect to the $L^2$-norm with constant one. Moreover, for $\pih$ we get by Poincaré's inequality
\cite{PaW60}
that, for all $K \in \Th$ and for all $v \in H^1(K)$,
\begin{equation}
	\label{eq:pwpoincare}
	\tnf{v-\pih v}{K} \leq \pi^{-1} h \tnf{\nabla v}{K},
\end{equation}
where we write $L^2(K)$ for the restriction of the $L^2$-space and its associated inner product and norm to the element~$K$. If no subdomain is specified, we always refer to the $L^2$-space on the whole domain. The following lemma from \cite{Gal23} establishes a crucial commuting property for the operators~$\pih$ and~$\Pih$.
\begin{lemma}[Commuting property]\label{lem:commut}
    Any $ v \in H^1_0(\Omega)$ satisfies
    $
		G_h \pih  v = \Pih\nabla v.
	$
\end{lemma}
\begin{proof}
    Using~\cref{eq:discgrad} and integration by parts, we obtain that 
	\begin{align*}
		\tspof{G_h \pih  v}{\tau_h}
		=
		-\tspof{\ddiv\tau_h}{\pih v}
		&=
		-\tspof{\ddiv\tau_h}{v}
		 =
		\tspof{ \nabla v}{\tau_h}
		=
		\tspof{\Pih\nabla v}{\tau_h}
	\end{align*}
for any $ v \in H^1_0(\Omega)$ and $\tau_h\in\Sigma_h$, which is the assertion.
\end{proof}

\section{Guaranteed lower energy bound}

The following theorem gives a lower bound on the ground state energy $E:=\mathcal{E}(u)$ using a post-processed version of the discrete ground state energy $E_h \coloneqq \mathcal E_h(u_h)$. This is the first major result of this paper.

\begin{theorem}[Lower bound]\label{th:lowerbound}
If the potential $V$ is $\Th$-piecewise constant, it holds that
	\begin{equation}
		\label{eq:lowerbound}
		\frac{E_h }
		{1+4 h^2 \pi^{-2}E_h }
		\leq E.
	\end{equation}
\end{theorem}

\begin{proof}
The discrete energy of the ground state is characterized by the following pseudo-Rayleigh quotient
	\begin{equation*}
		E_h
		=
		\min_{v_h\in U_h\setminus\{ 0\}}
		\frac{
			\tfrac12\tnof{G_h v_h}^2 \tnof{v_h}^2
			+  \tfrac12\tnof{V^{1/2} v_h}^2 \tnof{v_h}^2
			+ \tfrac{\kappa}{4}\|v_h\|_{L^4}^4
		}
		{\tnof{v_h}^4}.	
	\end{equation*}
We majorise $E_h$ by choosing $v_h:=\pih  u$. This results in 
\begin{equation}
	\label{eq:lowboundproof}
E_h \tnof{v_h}^4\leq
       \tfrac12 \tnof{G_h v_h}^2\tnof{v_h}^2
      +  \tfrac12\tnof{V^{1/2} v_h}^2 \tnof{v_h}^2
      + \tfrac{\kappa}{4} \|v_h\|_{L^4}^4.
\end{equation}
We bound all the terms on the right-hand side individually. Using the $L^2$-stability of~$\pih$, we get that
\begin{equation*}
\tnof{v_h}^2
=
\tnof{\pih  u}^2
\leq 
\tnof{u}^2.	
\end{equation*}
Since $V$ is assumed to be $\Th$-piecewise constant, we obtain that
\begin{equation*}
\tnof{V^{1/2} v_h}^2
=
\tnof{V^{1/2} \pih  u}^2
=
\sum_{K\in\Th}
V|_K \tn{\fint_K u \dx}{K}^2
\leq 
\tnof{V^{1/2} u}^2.	
\end{equation*}
\cref{lem:commut} and the $L^2$-stability of $\Pih$ yield that 
\begin{equation*}
	\tnof{G_h v_h}^2
	=
	\tnof{\Pih\nabla u}^2
	\leq 
	\tnof{\nabla u}^2.
\end{equation*}
Finally, the $L^4$-term is bounded by Jensen's
inequality
\begin{equation*}
	\|v_h\|_{L^4}^4
	=
	\sum_{K\in\Th}
	\int_K \bigg(\fint_K u \dx\bigg)^4   \dx
	\leq 
	\sum_{K\in\Th}
	\int_K \fint_K |u|^4 \dx  \dx
	=
	\|u\|_{L^4}^4.
\end{equation*}
Altogether, by inserting the above bounds into \cref{eq:lowboundproof} and using that $\tnof{u} = 1$, we get that
\begin{equation*}
	E_h \tnof{v_h}^4
	\leq E.
\end{equation*}
For rewriting the left-hand side we use the Pythagorean identity, which yields that
\begin{equation*}
\tnof{v_h}^4
=
\big(\tnof{\pih  u}^2 \big)^2
=
\big(1-\tnof{u-\pih  u}^2 \big)^2.
\end{equation*}
Using \cref{eq:pwpoincare}, we then obtain the lower bound
\begin{equation*}
\tnof{v_h}^4
\geq 
\big(1-h^2\pi^{-2}\tnof{\nabla u}^2 \big)^2\geq 
\big(1-2h^2\pi^{-2}E\big)^2,	
\end{equation*}
where we have used that  $\tfrac12\tnof{\nabla u}^2\leq E$. The combination of the previous estimates leads to the inequality
\begin{equation*}
	E_h \big(1-2h^2\pi^{-2}E\big)^2 \leq E.
\end{equation*}
Expanding the squared brackets and estimating yields
\begin{equation*}
	\big(1-2h^2\pi^{-2}E\big)^2 = 1
	-4h^2\pi^{-2}E
	+ 4h^4\pi^{-4}E^2
	\geq
	1
	-4h^2\pi^{-2}E,
\end{equation*}
which implies that  
\begin{equation}
	\label{eq:unifbounded}
	E_h 
	\big(1-4h^2\pi^{-2}E )
	\leq E.
\end{equation}
Elementary algebra then gives the assertion.
\end{proof}

\section{A priori error analysis}
In mixed form, the Gross-Pitaevskii eigenvalue problem for the ground state seeks the pair $(u,\sigma)\in L^2(\Omega)\times H(\operatorname{div},\Omega)$ with $\tnof{u} = 1$ corresponding to the smallest eigenvalue $\lambda \in \mathbb R$ such that
\begin{subequations}\label{e:o_evp}
	\begin{align}
			\tspof{\sigma}{\tau} + \tspof{\ddiv\tau}{u}
			&= 0 &&\text{for all }\tau\in H(\operatorname{div},\Omega),\label{e:o_evp_a}\\
			\tspof{\ddiv \sigma}{v} -  \tspof{(\kappa|u|^2+V) u}{v}
			&= -\lambda \tspof{u}{v} 
			&&\text{for all }v \in L^2(\Omega).\label{e:o_evp_b}
	\end{align}
\end{subequations}
Similarly, also any discrete ground state $u_h\in \Uh$ satisfies a mixed variational eigenvalue problem. More precisely, there exist $\sigma_h=G_hu_h\in \Sh$ and an eigenvalue $\lambda_h \in \mathbb R$ such that
\begin{subequations}\label{e:evp}
	\begin{align}
			\tspof{\sigma_h}{\tau_h} + \tspof{\ddiv\tau_h}{u_h}
			&= 0 &&\text{for all }\tau_h\in\Sigma_h,\label{eq:uhsigmah}\\
			\tspof{\ddiv \sigma_h}{v_h} -  \tspof{(\kappa|u_h|^2 +V)u_h}{v_h}
			&= -\lambda_h \tspof{u_h}{v_h} 
			&&\text{for all }v_h \in U_h.\label{eq:gpediscmixed}
	\end{align}
\end{subequations}
Note that $\lambda_h$ may not be the smallest discrete eigenvalue.  Similarly as in the continuous setting, the discrete energy and discrete ground state eigenvalue are related by $\lambda_h = 2E_h +\tfrac{\kappa}{2}\|u_h\|_{L^4}^4$.

The error analysis is based on the following elementary identity for the difference of the energies.
\begin{lemma}[Energy error characterization]\label{l:energydifference}
It holds that 
 \begin{equation*}
	\begin{aligned}
		E_h-E  
		&= -\tfrac12 \|G_h u_h -\nabla u\|_{L^2}^2
		    -\tfrac12 \| V^{1/2}(u_h-u)\|_{L^2}^2
		\\ 
		& \qquad
		    +\tfrac12 \lambda_h \|u_h-u\|_{L^2}^2
		    -\tfrac{\kappa}{4}\tspof{(u_h-u)^2}{3u_h^2+2uu_h+u^2}\\
		    &\qquad   + ((V-\pi_hV)(u_h-u),u_h)_{L^2}.
	\end{aligned}		
\end{equation*}
\end{lemma}
\begin{proof}
The definitions of $E_h$ and $E$ together with elementary
algebraic manipulations yield that
\begin{align*}
		E_h-E  
		&= -\tfrac12 \|G_h u_h -\nabla u\|_{L^2}^2
		    -\tfrac12 \| V^{1/2}(u_h-u)\|_{L^2}^2 + R
\end{align*}
with
\begin{align*}
R:=(G_h u_h-\nabla u, G_h u_h)_{L^2}
		    +(V(u_h-u),u_h)_{L^2}
		    +\tfrac{\kappa}{4} (\|u_h\|_{L^4}^4-\|u\|_{L^4}^4) .
\end{align*}
From the properties of the $L^2$-projections $\pih$ and~$\Pih$, the identity $\Pi_h\nabla u = G_h \pih u$ from \cref{lem:commut}, and~\cref{e:evp},
we get that
\begin{align*}
	 R&= (G_h u_h-\Pi_h\nabla u, G_h u_h)_{L^2}
	 +(\pi_hV(u_h-\pi_hu),u_h)_{L^2}\\* & \qquad+  ((V-\pi_hV)(u_h-u),u_h)_{L^2}
	 +\tfrac{\kappa}{4} (\|u_h\|_{L^4}^4-\|u\|_{L^4}^4)\\
	 &=  - \kappa(u_h^3,u_h-u)_{L^2}
	+ \lambda_h(u_h,u_h-u)_{L^2} \\ & \qquad + ((V-\pi_hV)(u_h-u),u_h)_{L^2} + \tfrac{\kappa}{4} (\|u_h\|_{L^4}^4-\|u\|_{L^4}^4).
\end{align*}
Since $u_h$ and $u$ are $L^2$-normalised, we have
$\lambda_h(u_h,u_h-u)_{L^2} =\tfrac12 \lambda_h \|u_h-u\|_{L^2}^2$. Rearranging the terms and using that
$$
- \kappa(u_h^3,u_h-u)_{L^2} + \tfrac{\kappa}{4} (\|u_h\|_{L^4}^4-\|u\|_{L^4}^4)
=
-\tfrac{\kappa}{4}\tspof{(u_h-u)^2}{3u_h^2+2uu_h+u^2},
$$
readily yields the assertion. 
\end{proof}

\begin{remark}[Tilde notation]
In the following, we will write $a \lesssim b$ or $b\gtrsim a$ if it holds that $a \leq C b$ or $a \geq C b$, respectively, where $C>0$ is a constant that may depend on the domain, the mesh regularity, the coefficients $V$ and $\kappa$, and on the ground state $u$, but is independent of the mesh size $h$. 
\end{remark}

The following theorem states a convergence result for the mixed finite element approximation to the ground state. 

\begin{theorem}[Plain convergence of mixed method]\label{thm:plainconv}
As $h\to 0$, it holds that
\begin{equation*}
	\tnof{u-u_h} \rightarrow 0,\quad \tnof{G_hu_h- \nabla u} \rightarrow 0,\quad E_h \rightarrow E,\quad \lambda_h\rightarrow \lambda.
\end{equation*}
\end{theorem}
\begin{proof}
We consider
\begin{equation}
	\label{eq:modi:energy}
	u_h^* \in  \argmin_{v \in H^1_0(\Omega)\with \tnof{v} = 1} \mathcal E_h^{*} (v)
\end{equation}
with the modified energy
\begin{equation*}
	\mathcal E_h^{*} (v)  \coloneqq \tfrac12\tspof{\nabla v}{\nabla v} + \tfrac12 \tspof{\pi_hVv}{v} + \tfrac{\kappa}{4} \tspof{|v|^2 v}{v},\quad v \in H^1_0(\Omega).
\end{equation*}
Similar as for the Gross-Pitaevskii energy minimisation problem~\cref{eq:gs}, the global modified energy minimiser exists and is unique up to sign. Note that despite the use of $h$ in the notation, \cref{eq:modi:energy} is a continuous problem. To get the uniqueness of \cref{eq:modi:energy}, we choose the sign of $u_h^*$ such that $(u,u_h^*)_{L^2}\geq 0$ holds.
The energies $E_h^*\coloneqq \mathcal E_h^*(u_h^*)$ are uniformly bounded with respect to $h$ since
\begin{equation*}
\label{eq:unifboundEhs}
	E_h^* \leq\mathcal E_h^*(u)=E + \tfrac{1}{2}\tspof{(\pi_hV-V)u}{u} \lesssim 1 + \tnof{V-\pih V}\lesssim 1,
\end{equation*}
where we used that $E_h^* \leq \mathcal E_h^*(u)$ and the $L^4$-regularity of $u$.  The uniform boundedness of $E_h^*$, directly implies that $\|u_h^*\|_{L^4}$ is uniformly bounded. Using this, we obtain similarly as before that 
\begin{equation}\label{eq:modi:approxi1}
	\begin{aligned}
			E_h^{*} - \mathcal E (u_h^*) = &\tfrac{1}{2}((\pi_hV-V)u_h^*,u_h^*)_{L^2}\rightarrow 0,\\  \mathcal E_h^{*} (u) - E = &\tfrac{1}{2}((\pi_hV-V)u,u)_{L^2}\rightarrow 0,
	\end{aligned}
\end{equation}
which together with $E \leq \mathcal E(u_h^*)$ and $E_h^* \leq \mathcal E_h^*(u)$ implies that 
\begin{equation}\label{eq:modi:approxi2}
	\begin{aligned}
			0 &\leq \mathcal E (u_h^*) - E = \mathcal E (u_h^*)- E_{h}^{*}  + E_h^{*}- \mathcal E^{*}_h (u) + \mathcal E^{*}_h (u)- E \\ &\leq |E_h^{*}- \mathcal E (u_h^*)| + |\mathcal E_h^{*} (u) - E| \to 0.
	\end{aligned}
\end{equation}
Combining \cref{eq:modi:approxi1,eq:modi:approxi2}, we get that
\begin{align}\label{eq:EhsmE}
	|E-E_h^*| \leq |E_h^*-\mathcal E (u_h^*)| + |\mathcal E (u_h^*) - E| \to 0.
\end{align}

Note that the discrete ground state can be interpreted as a discretisation of~\cref{eq:modi:energy}. This allows us to conclude, similarly to \cref{eq:unifbounded} in the proof of \cref{th:lowerbound}, that 
\begin{equation}
	\label{eq:modi:unifbounded}
	E_h 
	\big(1-4h^2\pi^{-2}E_h^*)
	\leq E_h^*,
\end{equation}
which implies the uniform boundedness of the discrete energies. As a consequence $\tnof{G_hu_h}$, $\|u_h\|_{L^4}$, and $\lambda_h$ are uniformly bounded. Furthermore, by the discrete embedding of \cref{l:dem}, $\|u_h\|_{L^6}$ is also uniformly bounded. 

Using the uniform bounds from above, we have that  
\begin{equation*}
	\tnof{\ddiv G_h u_h} = \tnof{\pi_h(\kappa|u_h|^2u_h + Vu_h - \lambda_h u_h)} \lesssim 1.
\end{equation*}
This estimate has two consequences: First, by \cref{lem:uh_Linfty}, it implies the uniform boundedness of $\|u_h\|_{L^\infty}$. Second, denoting by $\tilde u_h^c\in H^2(\Omega)\cap H_0^1(\Omega)$ the conforming lifting of $u_h$ from \cref{lem:appro:conti}, we have that 
\begin{equation}\label{eq:approx_plainproof}
	\tnof{G_hu_h-\nabla  \tilde u_h^c} + \tnof{u_h- \tilde u_h^c}\lesssim h  \tnof{\ddiv  G_hu_h} \lesssim h,
\end{equation}
where we used the bound from \cref{lem:appro:conti}. In the following, we consider the $L^2$-normalised version  $u_h^c  \coloneqq \tilde u_h^c/\|\tilde u_h^c\|_{L^2}$. Using elementary algebra one can show for the normalisation constant that 
\begin{equation}\label{eq:ele:algebra}
	\bigl|\|  \tilde u_h^c\|_{L^2} - 1 \bigr|
	= \bigl|\|  \tilde u_h^c\|_{L^2} - \|u_h\|_{L^2} \bigr|\leq \tnof{u_h- \tilde u_h^c} \lesssim h,
\end{equation}
where we used that $\|u_h\|_{L^2} =1$. 

Combining \cref{eq:ele:algebra,eq:approx_plainproof} and using  \cref{eq:L6bound}, the estimate 
\begin{equation}\label{eq:approxtilde}
	\tnof{G_hu_h-\nabla  u_h^c} + \tnof{u_h- u_h^c}\lesssim h
\end{equation}
is easily derived for sufficiently small $h>0$. Therefore $\|\nabla u_h^c\|_{L^2}$ and $\|u_h^c\|_{L^2}$ are uniformly bounded. By the embedding $H^1(\Omega)\hookrightarrow L^6(\Omega)$, $\|u_h^c\|_{L^6}$ and $\|u_h^c\|_{L^4}$ are also uniformly bounded. This implies that 
\begin{equation}\label{eq:E-mE-hstuhc}
	\begin{aligned}
		E_h - \mathcal E_h^{*}(u_h^c) &= \tfrac{1}{2}(G_hu_h-\nabla u_h^c,G_hu_h+\nabla u_h^c)_{L^2} + \tfrac{1}{2}(\pi_hV(u_h-u_h^c),u_h+u_h^c)_{L^2}\\ &\qquad + \tfrac{\kappa}{4}(\|u_h\|_{L^4}^4-\|u_h^c\|_{L^4}^4)\rightarrow 0 .
	\end{aligned}
\end{equation}
and 
\begin{equation}\label{eq:E-mE-hstuhc1}
		\mathcal E_h^{*}(u_h^c) - \mathcal E(u_h^c) = \tfrac{1}{2}((\pi_hV-V)u_h^c,u_h^c)_{L^2}\rightarrow 0.
\end{equation}
The inequality $E_h^{*}\leq \mathcal E_h^{*}(u_h^c)$, the lower bound \cref{eq:modi:unifbounded},  and \cref{eq:E-mE-hstuhc} imply that 
\begin{equation*}
	0\geq   E_h^{*} - \mathcal E_h^{*}(u_h^c)\geq E_h \big(1-4h^2\pi^{-2}E_h^*)- \mathcal E_h^{*}(u_h^c)\rightarrow 0.
\end{equation*} 
This together with~\cref{eq:EhsmE,eq:E-mE-hstuhc1} gives that
\begin{equation}\label{eq:estE-Etuhc}
	|E-\mathcal E(u_h^c)| \leq  |E- E_h^*| + |E_h^* - \mathcal E_h^{*}(u_h^c)|  + |\mathcal E_h^{*}(u_h^c) - \mathcal E(u_h^c)| \to 0.
\end{equation}
Assuming that $(u,u_h^c)_{L^2}\geq 0$ holds for $h$ sufficiently small, one can show that $\|u-u_h^c\|_{H^1}\rightarrow0$ using \cref{eq:estE-Etuhc} and similar arguments as in the proof of \cite[Thm.~1]{CCM10}. 
Otherwise, one can proceed with $v_h^c \coloneqq -u_h^c$, which similarly yields that $\|u-v_h^c\|_{H^1}\to 0$. Since on the one hand $\| u+u_h\|_{L^2}\leq \|u -v_h^c\|_{L^2} + \|u_h^c-u_h\|_{L^2}\to 0$ and on the other hand $\|u+u_h\|_{L^2}^2 = 2+ 2(u,u_h)_{L^2} \geq 2$, we get a contradiction which shows that it must hold that $(u,u_h^c)_{L^2}\geq 0$ for $h$ sufficiently small.

The  convergence  of the energies, i.e., $E_h \rightarrow E$, follows immediately combining \cref{eq:E-mE-hstuhc,eq:E-mE-hstuhc1,eq:estE-Etuhc}. To show the $L^2$-convergence of the gradient, we use the triangle inequality to obtain that
\begin{equation*}
	\tnof{G_hu_h-\nabla u} \leq \tnof{G_hu_h - \nabla u_h^c} + \tnof{\nabla u_h^c - \nabla u} \rightarrow 0.
\end{equation*}
Similarly, one can show that $\tnof{u - u_h}\rightarrow 0$. For the eigenvalues, we get that
\begin{equation}
	\label{eq:estlambda}
	|\lambda_h-\lambda| \leq 2|E_h-E| + \tfrac{\kappa}{2} |\|u_h\|_{L^4}^4 - \|u\|_{L^4}^4|\rightarrow 0.
\end{equation}
Algebraic manipulations and the application of Hölder's inequality for the second term on the right-hand side  prove the convergence of the ground state eigenvalue approximation. This concludes the proof. 
\end{proof}

The following corollary is an immediate consequence of the previous proof. 
\begin{corollary}[Uniform boundedness]\label{cor:unifbound}
	It holds that $\tnof{G_hu_h}$, $\|u_h\|_{L^\infty}$, $\lambda_h$, and~$E_h$ are uniformly bounded with respect to $h$.
\end{corollary}

For the quantification of the rates of convergence, we introduce some new notation. We denote pairs of functions in $L^2(\Omega)\times H(\ddiv,\Omega)$ by boldface
Roman capital letters, e.g. $\mathbf U$, $\mathbf V$, and $\mathbf W$.
The discrete analogues in $\Uh\times \Sh$ are denoted by $\mathbf U_h$, $\mathbf V_h$, and $\mathbf W_h$. 
Furthermore, we define the bilinear form $\bu$ acting on the pairs $\mathbf V = (v,\tau)$ and $\mathbf W = (w,\vartheta)$ as follows
\begin{equation*}
	\bu(\mathbf V,\mathbf W) \coloneqq \tspof{\tau}{\vartheta} + \tspof{\ddiv\vartheta}{v}-\tspof{\ddiv\tau}{w}+ \tspof{Vv}{w}+ \kappa\tspof{|u|^2v}{w},
\end{equation*}
where $u$ denotes the ground state. 
By rewriting~\cref{eq:gpe} as a Poisson problem with the $L^2$-right-hand side $\lambda u - Vu - \kappa|u|^2u$ and using the embedding $H^1(\Omega)\hookrightarrow L^6(\Omega)$ for $d \leq 3$, classical elliptic regularity theory (see, e.g., \cite[Thm.~9.1.22]{Hac03}) easily shows that the ground state~$u$ is $H^2$-regular, i.e., $u \in H^2(\Omega)\cap H^1_0(\Omega)$. The embedding $H^2(\Omega)\hookrightarrow C^0(\overline\Omega)$ for $d \leq 3$ then shows that $u$ is essentially bounded, i.e., its $L^\infty$-norm is finite. 
This in turn shows that the bilinear form $B_u$ is well defined. 

Similarly, denoting the ground state eigenvalue by $\lambda$, we define the bilinear form $\jul$ by
\begin{equation*}
	\jul(\mathbf V,\mathbf W)\coloneqq \bu(\mathbf V,\mathbf W) - \lambda \tspof{v}{w} + 2\kappa \tspof{|u|^2v}{w}.
\end{equation*}

We can then prove the following preliminary result.
\begin{lemma}[Almost coercivity of $\jul$]\label{lem:estimate:elliptic}
	For any $\mathbf V_h = (v_h,G_h v_h)$, it holds that
	\begin{equation*}
		\tnof{G_hv_h}^2 + \tnof{v_h}^2 \lesssim  \jul(\mathbf V_h,\mathbf V_h)+ h^2\tnof{\ddiv G_hv_h}^2.
	\end{equation*}
\end{lemma}
\begin{proof}
	By \cref{lem:appro:conti}, there exists for any $\mathbf V_h = (v_h,G_hv_h)$ a pair $\mathbf V_h^c = (v_h^c,\nabla v_h^c)$ with $v_h^c\in H^2(\Omega)\cap H_0^1(\Omega)$ such that it holds
	\begin{equation}
		\label{eq:approxconfcounterpart}
		\tnof{G_hv_h-\nabla v_h^c} + \tnof{v_h-v_h^c}\lesssim h  \tnof{\ddiv G_hv_h}.
	\end{equation}
Using \cite[Lem.~1]{CCM10}, which provides a lower bound of $\jul$ for conforming functions, we obtain that 
\begin{align*}
			\tnof{G_hv_h}^2 + \tnof{v_h}^2 &\lesssim   \tnof{\nabla v_h^c}^2 + \tnof{v_h^c}^2 + h^2\tnof{\ddiv G_hv_h}^2\\
			&\lesssim  \jul (\mathbf V_h^c,\mathbf V_h^c) +  h^2\tnof{\ddiv  G_hv_h}^2.\notag
\end{align*}
The desired result does not include the conforming counterpart $\mathbf V_h^c$ but $\mathbf V_h$. To go back to the original function $\mathbf V_h$, we use~\cref{eq:approxconfcounterpart}, the $L^\infty$-bound for $u$, and Young's inequality to get that 
	\begin{align*}
		 &| \jul (\mathbf V_h^c,\mathbf V_h^c) - \jul(\mathbf V_h,\mathbf V_h)| \\
		 &\qquad \leq|\tspof{\nabla v_h^c}{\nabla v_h^c - G_hv_h}| + |\tspof{G_hv_h}{\nabla v_h^c - G_hv_h}| \\
		  &\qquad\qquad   + |\tspof{(V+3\kappa u^2-\lambda)(v_h^c-v_h)}{v_h^c+v_h}|\\
		 &\qquad  \lesssim  h \tnof{\ddiv G_h v_h}\big(h\tnof{\ddiv G_h v_h} + \tnof{G_hv_h} + \tnof{v_h}\big)\\
		 &\qquad \leq h^2\big(1+\tfrac{1}{4\epsilon}\big) \tnof{\ddiv G_h v_h}^2 + \epsilon\big(\tnof{G_hv_h}^2 + \tnof{v_h}^2\big),
	\end{align*}
which holds for all $\epsilon>0$. 
Combining the previous two estimates yields
\begin{align*}
	&\tnof{G_hv_h}^2 + \tnof{v_h}^2\\
	&\qquad  \leq C\big( \jul (\mathbf V_h,\mathbf V_h) + h^2\big(1+\tfrac{1}{4\epsilon}\big) \tnof{\ddiv G_h v_h}^2 + \epsilon\big(\tnof{G_hv_h}^2 + \tnof{v_h}^2\big)\big)
\end{align*}
for some $C>0$ independent of $h$ and $v_h$. By choosing $\epsilon =\tfrac{1}{2C}$, the rightmost term can be absorbed into the left-hand side. This completes the proof.
\end{proof}

Let $\mathbf U = (u,\sigma)$ and $\mathbf U_h = (u_h,\sigma_h)$ denote solutions of
\eqref{e:o_evp} and \eqref{e:evp}, respectively. Recall that we assume $(u,u_h)_{L^2}\geq 0$ so that $u$ and $u_h$ have compatible signs. 
To simplify the notation, we introduce the $L^2$-norm 
in the product space for any ${\mathbf V = (v,\tau)}$ as 
$\|\mathbf V\|^2 \coloneqq \tnof{v}^2 + \tnof{\tau}^2$.
Let us define the pair $\tilde{\mathbf W}_h = (\tilde{w}_h,\tilde{\vartheta}_h)$
as the solution to 
\begin{subequations}
\begin{align}
	\tspof{\tilde{\vartheta}_h}{\tau_h} + \tspof{\ddiv\tau_h}{\tilde{w}_h} &=0 && \text{for all }\tau_h\in \Sh,\label{eq:tildew}\\*
	\tspof{\ddiv \tilde{\vartheta}_h}{v_h}  &= \tspof{\Delta u}{v_h} && \text{for all }v_h\in \Uh
\end{align}
\end{subequations} 
and set ${\mathbf W}_h = ({w}_h,{\vartheta}_h) = \tilde{\mathbf W}_h/\|\tilde{w}_h\|_{L^2}$. 
To prove an error estimate for $\|\mathbf U-\mathbf U_h\|$, we use the triangle inequality and examine the two errors 
$\|\mathbf U-\mathbf W_h\|$ and $\|\mathbf U_h-\mathbf W_h\|$ individually. 

\begin{lemma}
	[Estimate of first term]
	\label{lemma:errestI}
	For $h>0$ sufficiently small, it holds that
	\begin{equation*}
		\|\mathbf U-\mathbf W_h\|\lesssim h.
\end{equation*}
\end{lemma}
\begin{proof}
	A standard a~priori error estimate, cf.~\cite[Prop.~7.1.2]{BoffiBrezziFortin2013}, shows that
	\begin{equation*}
		\tnof{\tilde{\vartheta}_h-\nabla u} + \tnof{\tilde{w}_h-u}
		\lesssim h \|u\|_{H^2}
	\end{equation*}
since  $(\tilde{w}_h,\tilde{\vartheta}_h)$
	is the mixed Galerkin projection of $(u,\nabla u)$. The desired estimate then immediately follows from $\|u\|_{L^2} = 1$ and $|\|\tilde{w}_h\|_{L^2}-1| \lesssim h$.
\end{proof}

The following lemma is the final step towards the desired error estimate.

\begin{lemma}[Estimate of second term]\label{lemma:errestII}
	For $h>0$ sufficiently small, we have that
	\begin{equation*}
\|\mathbf U_h-\mathbf W_h\|\lesssim \|(u_h-u)^2\|_{L^3} + h\tnof{u_h-u} + h.
	\end{equation*}
\end{lemma}
\begin{proof}
We abbreviate $\mathbf Y_h:=\mathbf U_h - \mathbf W_h$ and $y_h:=u_h - w_h$.
	Using \cref{lem:estimate:elliptic}, we obtain that
	\begin{align*}
		\|\mathbf U_h - \mathbf W_h\|^2 
		&\lesssim \jul(\mathbf U_h-\mathbf W_h,\mathbf U_h-\mathbf W_h) 
		+ h^2\tnof{\ddiv (\sigma_h-\vartheta_h)}^2= \Xi_1 + \Xi_2+\Xi_3,
    \end{align*}
    where we set
    $$
            \Xi_1:=\jul(\mathbf U_h-\mathbf U,\mathbf Y_h),
            \quad
            \Xi_2:=\jul(\mathbf U-\mathbf W_h,\mathbf Y_h),
            \quad
            \Xi_3:= h^2\tnof{\ddiv (\sigma_h-\vartheta_h)}^2.
	$$
The term $\Xi_1$ is rewritten as follows
\begin{align*}
\Xi_1 
	&= \bu(\mathbf U_h-\mathbf U,\mathbf Y_h) - \lambda \tspof{u_h-u}{y_h} + 2\kappa\tspof{u^2(u_h - u)}{y_h}\\
	&=\lambda_h(u_h,y_h)_{L^2} - \kappa\tspof{(u_h^2-u^2)u_h}{y_h} - \lambda(u,y_h)_{L^2}\\ &\qquad \qquad- \lambda \tspof{u_h-u}{y_h} + 2\kappa\tspof{u^2(u_h - u)}{y_h}\\
	& =\tfrac12(\lambda_h-\lambda)\tnof{y_h}^2 - \kappa\tspof{(u_h^2-u^2)u_h}{y_h} + 2\kappa\tspof{u^2(u_h-u)}{
	y_h}\\
	& =\tfrac12(\lambda_h-\lambda)\tnof{y_h}^2 - \kappa\tspof{(u_h-u)^2(u_h+2u)}{y_h},
\end{align*}
where we use that $(u_h,y_h)_{L^2} = (u_h,u_h-w_h)_{L^2} = \tfrac{1}{2}\|u_h-w_h\|_{L^2}^2$, and  that $(u,\sigma)$ and $(u_h,\sigma_h)$ solve \eqref{e:o_evp} and \eqref{e:evp}, respectively. Using that  $\|u_h+2u\|_{L^6}\lesssim \tnof{G_hu_h} + \tnof{\nabla u}\lesssim 1$ (cf.~\cref{l:dem}), we then obtain the following estimate for $\Xi_1$
\begin{equation*}
	|\Xi_1|\lesssim |\lambda_h-\lambda|\tnof{y_h}^2 + \|(u_h-u)^2\|_{L^3}\tnof{y_h}.
\end{equation*}

For the term $\Xi_2$, we get that
\begin{align*}
	\Xi_2&= \bu(\mathbf U-\mathbf W_h,\mathbf Y_h) - \lambda \tspof{u-w_h}{y_h} + 2\kappa\tspof{u^2(u-w_h)}{y_h}\\
	& = (\sigma-\vartheta_h,\sigma_h-\vartheta_h)_{L^2} + \tspof{\ddiv (\sigma_h-\vartheta_h)}{u-w_h} -\tspof{\ddiv (\sigma- \vartheta_h)}{y_h}\\&\qquad\qquad\qquad + \tspof{(V+3\kappa u^2-\lambda)(u-w_h)}{y_h}\\
	& = \tspof{\ddiv (\vartheta_h- \tilde\vartheta_h)}{y_h}+ 
	\tspof{(V+3\kappa u^2-\lambda)(u-w_h)}{y_h}\\
	&=- \big(1-\tfrac{1}{\|\tilde w_h\|_{L^2}}\big)\tspof{(V+\kappa u^2-\lambda)u}{y_h} + 
	\tspof{(V+3\kappa u^2-\lambda)(u-w_h)}{y_h},
\end{align*}
where we used that 
$\tspof{\ddiv (\sigma-\tilde \vartheta_h)}{y_h}=0$ and $\Delta u = \kappa u^3+ Vu-\lambda u$, as well as
the identity $\tspof{\sigma-\vartheta_h}{\sigma_h-\vartheta_h} + \tspof{\ddiv (\sigma_h-\vartheta_h)}{u-w_h} = 0$, which is derived by integrating by parts and using~\cref{eq:tildew}. 
The estimate  $|\|\tilde{w}_h\|_{L^2}-1| \lesssim h$ then allows us to bound $\Xi_2$ as follows
\begin{align*}
	|\Xi_2|\lesssim h \tnof{y_h} + \tnof{u-w_h}\tnof{y_h}.
\end{align*}

For the term $\Xi_3$, we note that
\begin{equation*}
	\ddiv \sigma_h = \pih \big(\kappa u_h^3+ Vu_h - \lambda_h u_h\big),
	\qquad 
    \ddiv \vartheta_h = \pih \big(\kappa u^3 +  Vu - \lambda u\big)/\|\tilde w_h\|_{L^2},
\end{equation*}
where $\Delta u = \kappa u^3 + Vu - \lambda u$. This gives us 
\begin{align*}
	\ddiv (\sigma_h - \vartheta_h)&= \pih \Big(\kappa(u_h-u)(u_h^2+u_hu+u^2)+V(u_h-u)-\lambda(u_h-u)\\ 
	&\qquad  + (\lambda-\lambda_h)u_h+ \big(1-\tfrac{1}{\|\tilde w_h\|_{L^2}}\big)\Delta u\Big).
\end{align*}
Using the $L^2$-stability of $\pih$, the (uniform) $L^{\infty}$-bounds for $u$ and $u_h$ (cf.~\cref{cor:unifbound}), $\tnof{u_h} = 1$, and that $|\|\tilde{w}_h\|_{L^2}-1| \lesssim h$ yields that
\begin{align*}
	\Xi_3 \lesssim h^2\big(\tnof{u_h-u}^2 + |\lambda-\lambda_h|^2 + h^2\big).
\end{align*}

Combining the above estimates for $\Xi_1, \Xi_2$, and $\Xi_3$, we obtain that
\begin{align*}
	\|\mathbf U_h-\mathbf W_h\|^2 
	&\lesssim 
	\big( \|(u_h-u)^2\|_{L^3}+ \tnof{u-w_h}+h\big) \tnof{y_h}\\
	    &\qquad +h^2\big(\tnof{u_h-u}^2 + |\lambda-\lambda_h|^2 + h^2\big),
\end{align*}
where we absorbed the term $|\lambda_h-\lambda|\tnof{y_h}^2$ into the left-hand side, which is possible for sufficiently small $h>0$; see \cref{thm:plainconv}.
Using \cref{lemma:errestI} and the weighted Young's inequality, we obtain that 
\begin{align*}
	\|\mathbf U_h-\mathbf W_h\| &\lesssim \|(u_h-u)^2\|_{L^3} + h
	+ h \tnof{u_h-u}+ h|\lambda-\lambda_h|+h^2.
\end{align*}
The assertion then follows from the uniform boundedness of $\lambda_h$ (see \cref{cor:unifbound}) and from the fact that $h^2\lesssim h$  for $h>0$ sufficiently small.
\end{proof}

The following theorem gives an error estimate for the ground state, energy, and eigenvalue approximations of the proposed mixed finite element discretisation. It is derived by combining the two previous lemmas. For a second-order estimate for the eigenvalue approximation, which holds under additional regularity assumptions on $V$, we refer to  \cref{th:errestimp}.

\begin{theorem}[A priori error estimates]
	\label{th:errest}
	For sufficiently small $h>0$, it holds that
	\begin{equation}
		\label{eq:errest}
		\|\mathbf U-\mathbf U_h\|\lesssim h,\qquad \tnof{\ddiv (\sigma-\sigma_h)}\lesssim h + \|Vu-\pih (Vu)\|_{L^2}.
	\end{equation}
The eigenvalue and energy approximations satisfy
\begin{equation}
	\label{eq:estle} 
	|E-E_h| \lesssim h^2 + h\|V-\pih V\|_{L^2},\qquad |\lambda-\lambda_h| \lesssim h.
\end{equation}
If, in addition, $V$ is $\Th$-piecewise constant or $H^1$-regular, we have that
\begin{equation}
	\label{eq:impdivenergy}
	\tnof{\ddiv (\sigma-\sigma_h)}\lesssim h,\qquad |E-E_h| \lesssim h^2.
\end{equation}
\end{theorem}

\begin{proof}
	By the triangle inequality and \cref{lemma:errestI,lemma:errestII}, we obtain that
	\begin{align*}
		\|\mathbf U-\mathbf U_h\| 
		&\leq \|\mathbf U-\mathbf W_h\| + \|\mathbf U_h-\mathbf W_h\|\\
		&\lesssim \|(u_h-u)^2\|_{L^3} + h\tnof{u_h-u} + h.
	\end{align*}
For $h>0$ sufficiently small, the term $h\tnof{u_h-u}$ is absorbed into the left-hand side, which yields that
\begin{equation*}
	\|\mathbf U-\mathbf U_h\|\lesssim \|(u_h-u)^2\|_{L^3}+ h. 
\end{equation*}
It only remains to bound the first term on the right-hand side. Elementary algebraic manipulations and the triangle inequality show that
\begin{align}
	\label{eq:higherordertermest}
	\|(u_h-u)^2\|_{L^3} = \|u_h-u\|_{L^6}^2 \lesssim \|u_h - \pih u\|_{L^6}^2 + \|u-\pih u\|_{L^6}^2.
\end{align}
For the first term on the right-hand side of \cref{eq:higherordertermest}, we get with \cref{l:dem,lem:commut} and the triangle inequality that
\begin{align}
	\label{eq:l6normest}
	\begin{split}
			\|u_h - \pih u\|_{L^6}^2 &\lesssim \tnof{G_h (u_h - \pih u)}^2 \leq \tnof{G_h u_h - \nabla u}^2 + \tnof{\nabla u - \Pih \nabla u}^2\\
		&\lesssim \|\mathbf U-\mathbf U_h\|^2 + h^2,
	\end{split}
\end{align}
where we used the classical approximation property
\begin{equation*}
	\tnof{\nabla u - \Pih \nabla u} \lesssim h\|u\|_{H^2};
\end{equation*}
see, e.g.,~\cite[Prop.~2.5.4]{BoffiBrezziFortin2013}. The second term on the right-hand side of~\cref{eq:higherordertermest} is bounded using Poincaré's inequality \cite[Eq.~(7.45)]{GiT01} and applying the embedding $H^1(\Omega)\hookrightarrow L^6(\Omega)$ for $d \leq 3$ to the gradient of $u$. This results in
\begin{align}
	\label{eq:l6poincare}
	\|u-\pih u\|_{L^6}^2 \lesssim h^2 \|\nabla u\|_{L^6}^2 \lesssim h^2 \|u\|_{H^2}^2.
\end{align}
Combining the previous estimates, we obtain that
\begin{align*}
	\|\mathbf U-\mathbf U_h\|\lesssim \|\mathbf U-\mathbf U_h\|^2 + h^2+ h. 
\end{align*}
By the convergence result of \cref{thm:plainconv}, the term $\|\mathbf U-\mathbf U_h\|^2$ converges to zero. Since it is a higher order term, it is absorbed in the left-hand side for sufficiently small~$h>0$. Under this smallness condition, it also holds that $h^2 \lesssim h$. The desired estimate for $\|\mathbf U-\mathbf U_h\|$ follows immediately. 

The estimate for $|E-E_h|$ is an immediate consequence of \cref{l:energydifference}, the argument used in~\eqref{eq:higherordertermest}, and the uniform $L^\infty$-boundedness of $u_h$, cf.~\cref{cor:unifbound}. For proving the eigenvalue approximation result~\cref{eq:estle}, one may proceed similarly as in~\cref{eq:estlambda} using the convergence results for $\|\mathbf U-\mathbf U_h\|$ and $|E-E_h|$, cf.~\cref{eq:errest,eq:estle}.

Let us next prove the estimate for $\tnof{\ddiv (\sigma-\sigma_h)}$. We note that 
\begin{equation}
	\label{eq:dividentity}
	\ddiv \sigma_h = \pih \big(\kappa u_h^3+Vu_h - \lambda_h u_h\big),
	\qquad 
	\ddiv \sigma =  \kappa u^3 + Vu - \lambda u,
\end{equation}
which implies that
\begin{align*}
	\ddiv (\sigma_h-\sigma) & = \pi_h(\kappa(u_h^3-u^3)+V(u_h-u) -\lambda(u_h-u) -(\lambda_h-\lambda)u_h)\\ &\qquad + (\kappa(\pi_h u^3-u^3) + \pi_h(Vu) -Vu - \lambda(\pi_h u-u)).
\end{align*}
The desired estimate follows immediately using \cref{eq:pwpoincare}, the first estimate in \cref{eq:estle}, and \cref{eq:errest}. Finally, estimate \cref{eq:impdivenergy} is a direct consequence of \cref{eq:estle,eq:errest}.
\end{proof}

For $H^1$-regular potentials, the following theorem proves a second-order convergence result for the eigenvalue approximation.

\begin{theorem}[Improved error estimate]
	\label{th:errestimp}
	Let $V \in H^1(\Omega)$ and assume that $\Omega$ is a $d$-dimensional brick. Then, for sufficiently small $h>0$, it holds that
	\begin{equation}
		\label{eq:estleimp}
		\|u-u_h\|_{H^{-1}}\lesssim h^2,\qquad |\lambda-\lambda_h| \lesssim h^2.
	\end{equation}
\end{theorem}
\begin{proof}
	We begin with the proof of the $H^{-1}$-norm estimate.  For mixed finite elements, such error estimates were 
	introduced in \cite{DouglasRoberts1985}. In the following, we use this well-known technique with the auxiliary dual problem of 
	\cite[Eq.~(70)]{CCM10}.
	Given a test function $w\in H^1_0(\Omega)$, it seeks $z\in H^1_0(\Omega)$ such that 
	\begin{equation}
		\label{eq:adjprob}
		-\Delta z + (V+3\kappa u^2 -\lambda) z
		=
		2 \kappa(u^3,z)_{L^2} u + w - (w,u)_{L^2} u
	\end{equation}
	holds in $H^{-1}(\Omega)$. This problem is solved by the unique solution $z \in u^\perp \coloneqq \{v \in H^1_0(\Omega)\with \tspof{u}{v} = 0\}\subset H^1_0(\Omega)$ satisfying
	\begin{equation*}
		\tspof{\nabla z}{\nabla v} + \tspof{(V+3\kappa u^2-\lambda)z}{v} = \tspof{w}{v}\qquad \text{for all }v \in u^\perp.	
	\end{equation*}
	The well-posedness of the latter problem is a consequence of the Lax-Milgram theorem using the coercivity of the bilinear form on the left-hand side, cf.~\cite[Lem.~1]{CCM10}, and the fact that $u^\perp$ is a complete subspace of $H^1_0(\Omega)$. Assuming that $V \in H^1(\Omega)$, elliptic regularity theory implies that $z \in H^3(\Omega)$ with the estimate $\|z\|_{H^3} \lesssim \|w\|_{H^1}$. To prove the $H^3$-regularity, we recall the assumption that $\Omega$ is a $d$-dimensional brick and apply a prolongation by reflection argument, noting that the right-hand side of \cref{eq:adjprob} satisfies zero Dirichlet boundary conditions; see \cite[p.~107]{CCM10} for more details.

	Considering the mixed form of problem~\cref{eq:adjprob}, the pair $(z,\varphi)$ 
	satisfies
	\begin{subequations}\label{eq:mixedadjprob}
		\begin{align}
			\tspof{\varphi}{\tau} + \tspof{\ddiv\tau}{z}
			&= 0 && \text{for all }\tau\in H(\operatorname{div},\Omega),\label{eq:mixedadj_a}
			\\
			\tspof{\ddiv \varphi}{v} -  \tspof{(V+3\kappa u^2-\lambda)z}{v}&=\tspof{f}{v}\label{eq:mixedadj_b}
			&&\text{for all }v \in L^2(\Omega)
		\end{align}
	\end{subequations}
	for the source term
	$$
	f:=-2\kappa\tspof{u^3}{z} u - w +\tspof{w}{u} u.
	$$
	
	To derive the desired $H^{-1}$-norm estimate, we fix a test function $w \in H^1_0(\Omega)$ and test \cref{eq:mixedadj_b} with $u_h-u\in L^2(\Omega)$ and \cref{eq:mixedadj_a} with $\sigma_h-\sigma \in H(\ddiv,\Omega)$ and add up the equations. After rearranging the terms, we obtain that
	\begin{align*}
		\tspof{u-u_h}{w} &=  \Xi_1 + \Xi_2
	\end{align*}
	with the expressions
	\begin{align*}
		\Xi_1:=& \tspof{\varphi}{\sigma_h-\sigma} + \tspof{\ddiv (\sigma_h - \sigma)}{z} 
		+\tspof{\ddiv \varphi-(V+3\kappa u^2-\lambda)z}{u_h-u} 
		,
		\\
		\Xi_2:=&2\kappa\tspof{u^3}{z}\tspof{u}{u_h-u} - \tspof{w}{u}\tspof{u}{u_h-u}.
	\end{align*}
	In the following, we will add and subtract the term~$\tilde \Xi_1$ defined by
	\begin{align*}
		\tilde \Xi_{1} &\coloneqq  \tspof{\varphi_h}{\sigma_h-\sigma} + \tspof{\ddiv (\sigma_h - \sigma)}{z_h} 
		+\tspof{\ddiv \varphi_h-(V+3\kappa u^2-\lambda)z_h}{u_h-u}
	\end{align*}
	with $(z_h,\varphi_h) \coloneqq (\pih z,G_h \pih z)$. Let us first show that $|\tilde \Xi_1|$ is in fact a second-order term, i.e., $|\tilde\Xi_1|\lesssim h^2 \|w\|_{H^1}$. To see this, we seek a different representation of $\tilde \Xi_1$. Adding up~\cref{e:o_evp_a} tested with $\varphi_h$ and \cref{e:o_evp_b} tested with $z_h$ yields that
	\begin{align*}
		\tspof{\sigma}{\varphi_h} + \tspof{\ddiv \varphi_h}{u} + \tspof{\ddiv \sigma}{z_h} - \tspof{(V+\kappa u^2)u}{z_h} = -\lambda\tspof{u}{z_h}.
	\end{align*}
	Similarly, we get by adding up \cref{eq:uhsigmah} tested with  $\varphi_h$ and \cref{eq:gpediscmixed} tested with $z_h$ that
	\begin{align*}
		\tspof{\sigma_h}{\varphi_h} + \tspof{\ddiv \varphi_h}{u_h} + \tspof{\ddiv \sigma_h}{z_h} - \tspof{(V+\kappa u_h^2)u_h}{z_h} = -\lambda_h\tspof{u_h}{z_h}.
	\end{align*}
	Using these identities, we can rewrite $\tilde \Xi_1$ as
	\begin{align}
		\label{eq:tildexi}
		\tilde \Xi_1 = \kappa\tspof{(u_h^2+u_hu-2u^2)(u_h-u)}{z_h} + (\lambda-\lambda_h)\tspof{u_h}{z_h}.
	\end{align}
	Since it holds that $(u_h^2+u_hu-2u^2)(u_h-u) = (u_h+2u)(u_h-u)^2$, we obtain for the first term on the right-hand side of the previous equation  that
	\begin{align*}
		|\tspof{(u_h^2+u_hu-2u^2)(u_h-u)}{z_h}|\lesssim \|(u_h-u)^2\|_{L^3}\tnof{z_h}\lesssim 
		h^2\|w\|_{H^1},
	\end{align*}
	where we proceeded similarly as in~\cref{eq:higherordertermest} and used~\cref{l:dem}. For the second term on the right-hand side of~\cref{eq:tildexi}, we get with  $\tspof{z}{u} = 0$ that
	\begin{align*}
		(\lambda_h-\lambda)\tspof{u_h}{z_h} = (\lambda_h-\lambda)\tspof{u_h-u}{z_h} - (\lambda_h-\lambda)\tspof{u}{z-z_h},
	\end{align*}
	which, using \cref{eq:pwpoincare,eq:estle,eq:errest}, yields that
	\begin{align*}
		|(\lambda_h-\lambda)\tspof{u_h}{y_h}| \lesssim h^2\tnof{z_h} + h^2\|z\|_{H^1}
		\lesssim h^2\|w\|_{H^1}.
	\end{align*}
	
	Let us next estimate $|\Xi_1-\tilde \Xi_1|$. We use elementary algebraic manipulations to get that
	\begin{align*}
		\Xi_1-\tilde \Xi_1 &=\tspof{\varphi-\varphi_h}{\sigma_h-\sigma} + \tspof{\ddiv (\sigma_h-\sigma)}{z-z_h} + \tspof{\ddiv (\varphi-\varphi_h)}{u_h-u}\\
		& \qquad-  \tspof{(V+3\kappa u^2-\lambda)(u_h-u)}{z-z_h}.
	\end{align*}
	In the following, we estimate all terms on the right-hand side separately. For the first term, we get with \cref{lem:commut}, a classical approximation result, cf.~\cite[Prop.~2.5.4]{BoffiBrezziFortin2013}, and~\cref{eq:errest} that
	\begin{align*}
		|\tspof{\varphi-\varphi_h}{\sigma_h-\sigma}| \leq \tnof{\varphi - \Pih \varphi}\|\sigma-\sigma_h\|_{L^2} \lesssim h^2 \|w\|_{H^1}.
	\end{align*} 
	For the second term, a similar estimate can be obtained using \cref{eq:impdivenergy,eq:pwpoincare}.   
	
	Denoting by $I_h\colon H(\ddiv,\Omega)\to \Sigma_h$ the Raviart-Thomas interpolation operator, cf.~\cite[Sec.~2.5.2]{BoffiBrezziFortin2013}, we obtain for the third term that 
	\begin{align*}
		&\tspof{\ddiv (\varphi-\varphi_h)}{u_h-u} = \tspof{\ddiv (\varphi-I_h \varphi)}{u_h-u} + \tspof{\ddiv (I_h \varphi-\Pih \varphi)}{u_h-u}\\
		&\qquad =\tspof{\ddiv \varphi - \pih \ddiv \varphi}{u_h-u} - \tspof{\varphi - I_h \varphi}{\sigma-\sigma_h} + \tspof{\varphi-\Pih \varphi}{\sigma-\sigma_h}.
	\end{align*}
	Using classical approximation results for $I_h$ and $\Pih$, cf.~\cite[Prop.~2.5.4]{BoffiBrezziFortin2013}, and~\cref{eq:pwpoincare,eq:errest}, we obtain the estimate
	\begin{align*}
		|\tspof{\ddiv (\varphi-\varphi_h)}{u_h-u}| 
		\lesssim h\|\mathbf U-\mathbf U_h\|\|\varphi\|_{H^2}\lesssim h^2\|w\|_{H^1}.
	\end{align*}
	For the last term, we get using \cref{eq:pwpoincare,eq:errest} that
	\begin{align*}
		|\tspof{(V+3\kappa u^2-\lambda)(u_h-u)}{z-z_h}| \lesssim h\|u-u_h\|_{L^2}\|\nabla z\|_{L^2}\lesssim h^2\|w\|_{H^1}.
	\end{align*}
	Combining the previous estimates yields that
	\begin{equation*}
		|\Xi_1-\tilde \Xi_1| \lesssim h^2\|w\|_{H^1}.
	\end{equation*}
	
	For the term $\Xi_2$ we get with $\tspof{u}{u_h-u} = -\frac{1}{2} \|u_h-u\|_{L^2}^2$ the estimate 
	\begin{align*}
		|\Xi_2| \lesssim \big(\tnof{z} + \tnof{w}\big)\tnof{u_h-u}^2
		\lesssim h^2 \|w\|_{H^1}.
	\end{align*}
	Using that $|\tspof{u-u_h}{w}| \leq  |\tilde \Xi_1| + |\Xi_1-\tilde \Xi_1| + |\Xi_2| \lesssim h^2 \|w\|_{H^1}$ yields the desired estimate
	\begin{align*}
		\|u-u_h\|_{H^{-1}} = \sup_{w \in H^1_0(\Omega)\with \|w\|_{H^1} = 1}\tspof{u_h-u}{w} \lesssim h^2.
	\end{align*}
	
	Finally, to prove the second-order estimate for $|\lambda-\lambda_h|$, we  introduce the notation $c \coloneqq \kappa u^2+V$ and $c_h\coloneqq \kappa u_h^2+\pi_hV$. Using the identities
	\begin{align*}
		\tnof{\sigma}^2 &= \tnof{\sigma-\sigma_h}^2  + 2\tspof{\sigma}{\sigma_h} - \tnof{\sigma_h}^2,\\
		\tnof{c^{1/2} u}^2 &= \tnof{c^{1/2}(u-u_h)}^2 + 2\tspof{cu}{u_h} - \tnof{c^{1/2}u_h}^2
	\end{align*}
	and~\cref{e:o_evp_a,e:o_evp_b,eq:uhsigmah,eq:gpediscmixed},  algebraic manipulations yield that
	\begin{align*}
		\lambda - \lambda_h &= \tspof{-\ddiv \sigma + cu}{u} - \tspof{-\ddiv \sigma_h + c_h u_h}{u_h}\\
		&= \tnof{\sigma}^2 + \tnof{c^{1/2}u}^2 - \tnof{\sigma_h}^2 - \tnof{c_h^{1/2}u_h}^2\\
		&= \tnof{\sigma-\sigma_h}^2  + \tnof{c^{1/2}(u-u_h)}^2\\
		&\qquad + 2\tspof{\sigma - \sigma_h}{\sigma_h} 
		+ 2\tspof{cu}{u_h} - \tspof{c u_h}{u_h} - \tspof{c_hu_h}{u_h}.
	\end{align*}
	Using~\cref{e:o_evp_a,eq:uhsigmah,eq:gpediscmixed}, we get that
	\begin{align*}
		\tspof{\sigma-\sigma_h}{\sigma_h} &= \tspof{-\ddiv \sigma_h}{u-u_h} = \lambda_h \tspof{u_h}{\pi_hu-u_h} - \tspof{c_hu_h}{\pi_hu-u_h}\\
		&=-\tfrac{\lambda_h}{2}\tnof{u-u_h}^2 - \tspof{c_hu_h}{u-u_h},
	\end{align*}
	which yields the identity
	\begin{align*}
		\lambda-\lambda_h 
		&=\tnof{\sigma-\sigma_h}^2  + \tnof{c^{1/2}(u-u_h)}^2 - \lambda_h \tnof{u-u_h}^2 \\
		&\qquad + 2\tspof{(c-c_h)(u-u_h)}{u_h} +\tspof{(c-c_h)u_h}{u_h}. 
	\end{align*}
	Noting that $c-c_h =  \kappa(u-u_h)(u+u_h) + V-\pi_hV$ and using \cref{eq:pwpoincare} as well as the uniform $L^6$- and $L^\infty$-bounds for $u_h$ (cf.~\cref{cor:unifbound}), one obtains similarly as in~\cref{eq:higherordertermest} that 
	\begin{equation*}
		|\tspof{(c-c_h)(u-u_h)}{u_h}| \lesssim h^2.
	\end{equation*}
	Therefore, in order to show the second-order estimate for $|\lambda-\lambda_h|$, it only remains to consider the term
	\begin{equation*}
		\tspof{(c-c_h)u_h}{u_h} = \kappa\int_\Omega (u+u_h)(u-u_h)u_h^2\dx\eqqcolon \kappa\Xi.
	\end{equation*}
	Regularizing $u_h$ with the averaging operator $J$ from \cref{l:dem} yields that
	\begin{align*}
		\Xi &=  
		\int_\Omega (u-u_h) (u_h-J u_h) u_h^2 \dx+ \int_\Omega (u-u_h) (J u_h+u) (u_h^2-(Ju_h)^2) \dx\\
		&\qquad  +
		\int_\Omega (u-u_h) (J u_h+u) (Ju_h)^2 \dx.
	\end{align*}
	Noting that the gradient of $(J u_h+u) (Ju_h)^2$ can be computed as
	\begin{equation*}
		\nabla ( (J u_h+u) (Ju_h)^2 )
		=
		(3 J u_h^2 + 2uJu_h) \nabla Ju_h
		+
		(Ju_h)^2 \nabla u,
	\end{equation*}
	we obtain with 
	\cref{lem:properties_J} and the uniform $L^\infty$-bound of $u_h$ that
	\begin{equation*}
		\|(J u_h+u) (Ju_h)^2 \|_{H^1}\lesssim 1.	
	\end{equation*}
	With this, using the uniform $L^\infty$-bound of $u_h$ and $Ju_h$, as well as the approximation error estimate $\|u_h-Ju_h\|_{L^2}$ and $\|u-u_h\|_{L^2}$, and the first estimate in \cref{eq:estleimp}, we obtain that
	\begin{equation*}
		|\Xi|
		\lesssim
		h\|u-u_h\|_{L^2}
		+
		\|u-u_h\|_{H^{-1}}\lesssim h^2.
	\end{equation*}
	Combining the above estimates, the desired second-order approximation for $|\lambda-\lambda_h|$ immediately follows.
\end{proof}
\section{Numerical experiments}
Having laid the groundwork with our theoretical framework and error analysis for the mixed finite element discretisation of the Gross-Pitaevskii eigenvalue problem, we now shift our focus to numerical experiments. These experiments are essential both to validate our theoretical insights and to demonstrate the practicality of our approach. For the implementation, we have chosen solvers tailored to the finite-dimensional nonlinear eigenvector problem \cref{eq:gsdisc}, with the goal of aligning our numerical methods with the theoretical principles previously discussed.

In the realm of suitable methods, the discrete normalised gradient flow method referenced in \cite{BaD04} is a notable choice. This method is part of a diverse array of gradient flow approaches, each varying in their choice of metric, as indicated in \cite{RSS09,DaK10,KaE10}. An interesting advancement in this field is the introduction of an energy-adaptive metric, detailed in \cite{HenP20}, which has been further analysed for quantitative errors in subsequent studies \cite{Zha21ppt,CheLLZ23, AltmannPeterseimStykel2022}. Relatedly, Riemannian optimisation techniques, including Riemannian conjugate gradient \cite{DanP17,ALT17} and Riemannian Newton methods \cite{altmann2023riemannian}, offer additional avenues for exploration.
Other methods that focus on the formulation of the eigenvalue problem, such as the self-consistent field (SCF) iteration \cite{Can00,DioC07} and Newton's method \cite{JarU22}, also contribute valuable perspectives. It is interesting to note that assumptions about the symmetry of the condensate can lead to a reduction in the dimension of the problem, as explored in \cite{BaoT03}. Furthermore, the complexity of solving the nonlinear constraint minimisation problem can be reduced by using problem-adapted basis functions with high approximation quality \cite{HMP14b,henning2023optimal,PWZ23}, using techniques from (Super-)Localized Orthogonal Decomposition \cite{MaP14,HaPe21b}.

In this paper, we use the $J$-method of \cite{JarKM14,AltHP21} to solve the nonlinear discrete problem because, through the choice of shift, it nicely blends between the reliable linear convergence of gradient-descent type schemes and the local quadratic convergence of Newton-type methods. To apply the $J$-method in the mixed setting, we eliminate the dual variable in~\cref{e:evp}. This results in a system matrix of the form $(M(u) + CB^{-1}C^T)$, where $B$ is the Raviart-Thomas mass matrix, $C$ is the Raviart-Thomas divergence matrix, and $M(u)$ is a diagonal matrix containing the nonlinearity  and the potential. To avoid the costly computation of the Schur complement when solving with the system matrix, we use the Woodbury matrix identity. This gives
\begin{equation*}
	(M(u) + CB^{-1}C^T)^{-1} = M(u)^{-1} - M(u)^{-1}C(B+C^TM(u)^{-1}C)^{-1}C^TM(u)^{-1},
\end{equation*}
where the latter matrix is much easier to compute since 
$M(u)$ is diagonal. Note that since $M(u)$ is diagonal, $C^TM(u)^{-1}C$ is in fact a sparse matrix. For the damping, shifting, tolerances, etc., we use a similar parameter setting as in~\cite[Sec.~6]{AltHP21}. In particular, we use a damping strategy with an energy-diminishing step-size control when the $L^2$-norm of the residuals is larger than $10^{-2}$. For smaller residuals, damping is disabled and shifting is enabled. At this point, the method takes about three to four iterations to converge to machine accuracy. For implementation details, see the code provided at \url{https://github.com/moimmahauck/GPE_RT0}.

This section consists of two parts. First, we numerically investigate the optimal order convergence of the proposed mixed finite element discretisation of the Gross-Pitaevskii problem (see \cref{th:errest,th:errestimp}). Second, we numerically validate the lower bounds of the ground state energy (see \cref{th:lowerbound}).

\subsection*{Validation of optimal convergence rates}
\begin{figure}
	\includegraphics[height=.2\linewidth]{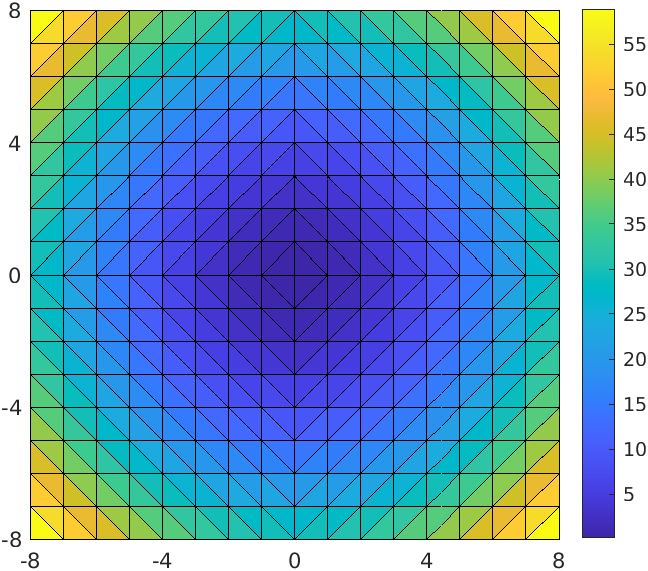}\hfill
	\includegraphics[height=.2\linewidth]{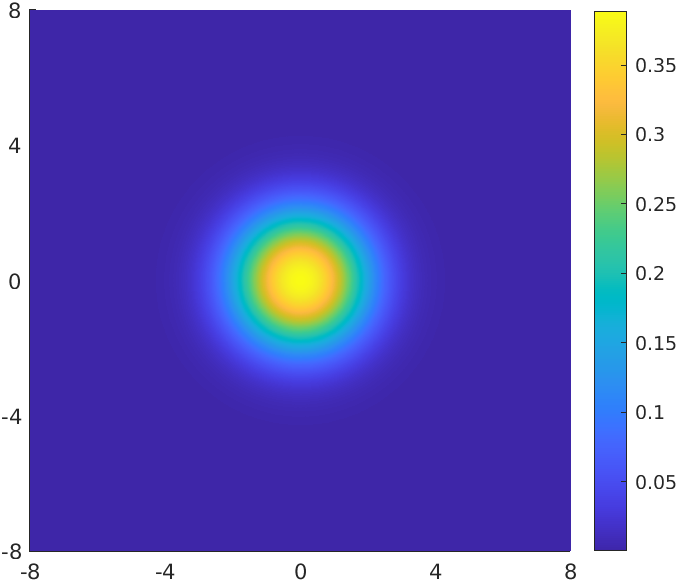}\hfill
	\includegraphics[height=.2\linewidth]{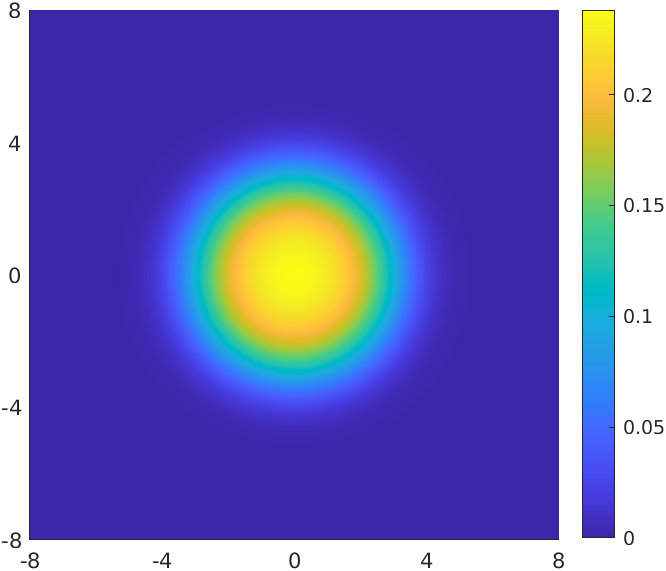}\hfill
	\includegraphics[height=.2\linewidth]{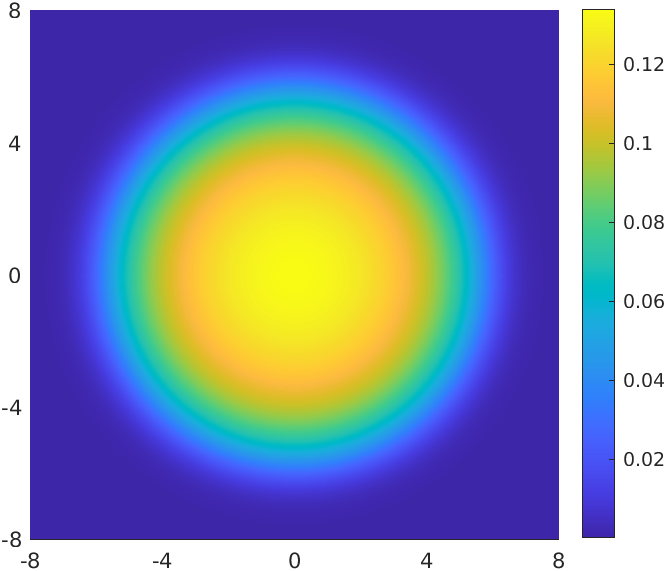}
	\caption{Projection of the potential onto the space of piecewise constants (left) and ground states for the harmonic potential for the parameters $\kappa = 10,100,1000$ (second to last plot).}
	\label{fig:groundstategarmpot}
\end{figure}

To verify the optimal order convergence, we consider the domain $\Omega = (-L,L)^2$ with $L = 8$ and the harmonic potential $V(x) = \tfrac12|x|^2$. 
For this setting, the ground state is point symmetric with respect to the origin and decays exponentially. The decay depends on the parameter $\kappa$: the larger~$\kappa$, the more repulsive the particle interaction and the more spread out the mass; see \cref{fig:groundstategarmpot} (last three plots). 
For the discretisation, we consider a hierarchy of meshes constructed by uniform red refinement of an initial mesh. The initial mesh is constructed from a Friedrichs-Keller triangulation consisting of eight elements by rotating the triangles in the lower right and upper left squares so that the mesh is point symmetric with respect to the origin. For each mesh in the hierarchy, we compute a ground state approximation, where we project the potential onto the space of piecewise constants with respect to the considered mesh; see \cref{fig:groundstategarmpot}~(left) for one projected potential. 

\begin{figure}
		\includegraphics[height=.45\linewidth]{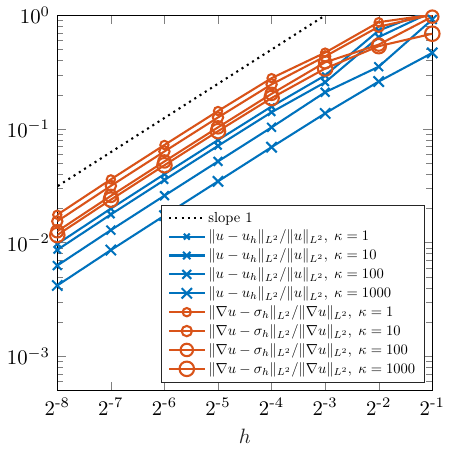}\hspace{1cm}
		\includegraphics[height=.45\linewidth]{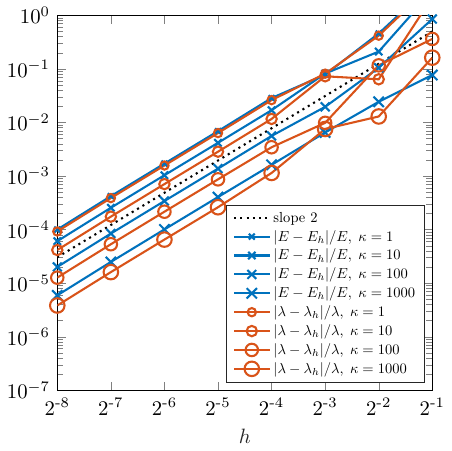}
		\caption{Error plots of the primal and dual variables (left) and of the energy and eigenvalue (right). The expected orders of convergence are indicated by black dotted lines.}
		\label{fig:convergencekappa}
\end{figure}

 \cref{fig:convergencekappa} then shows the errors of the mixed finite element discretisation for several values of $\kappa$. Note that since no analytical solution is available, the errors are computed with respect to a reference solution. This reference solution is computed on a mesh obtained by twice uniform red refinement of the finest mesh in the hierarchy. One observes first-order convergence for the primal and dual variables and second-order convergence for the energies and eigenvalues. Recalling that $V \in H^1(\Omega)$, this is consistent with the predictions in \cref{th:errest,th:errestimp}. We observe only a weak dependence of the errors on the parameter $\kappa$. More precisely, the errors are slightly smaller for larger~$\kappa$.

\subsection*{Validation of lower energy bounds}
Next, we numerically verify the lower ground state energy bound given in \cref {th:lowerbound}. We consider several different settings, namely a harmonic potential, a disorder potential, and a constant potential. In the following, we denote the post-processed discrete energy defined as the left-hand side of~\cref{eq:lowerbound} by $E_h^\mathrm{pp}$. To satisfy the assumption of \cref{th:lowerbound} that the potential is piecewise constant, we construct the potentials by prolongation of a piecewise constant potential on a coarse mesh. Reference values for the energies are computed using a $\mathcal Q^2$-finite element  implementation together with the energy-adaptive Riemannian gradient descent method from \cite{HenP20}. Note that, in order to use the same potentials for both methods, we choose the potential to be piecewise constant on a Cartesian mesh. For all our numerical experiments such a Cartesian mesh is constructed by joining opposing pairs of triangles of the coarse triangulation.

\subsubsection*{Harmonic potential with strong interaction}

First let us consider the harmonic potential $V(x) = \tfrac12|x|^2$ and the large parameter ${\kappa = 1000}$. The coarse mesh used for this numerical example is shown in the background of \cref{fig:groundstategarmpot} (left). 

\begin{figure}
	\includegraphics[width=.45\linewidth]{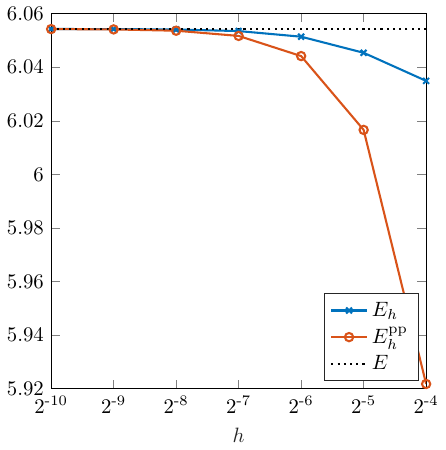}
	\hspace{.75cm}
	\includegraphics[width=.455\linewidth]{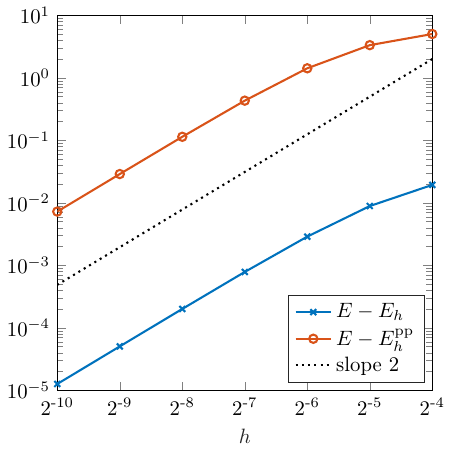}
	\caption{Ground state energy approximations for the harmonic potential (left) and the difference between the reference energy and the energy approximations in a double-logarithmic plot (right). The blue and red curves correspond to the discrete energy and the post-processed discrete energy, respectively.}
	\label{fig:lowboundharmpot}
\end{figure}

In \cref{fig:lowboundharmpot} (left) one observes that the energy $E_h$ and the post-processed energy $E_h^\mathrm{pp}$ strictly increase as $h$ is decreased, i.e., they approach the ground state energy from below. The observation  was predicted for the post-processed discrete energy by \cref{th:lowerbound}.  \cref{fig:lowboundharmpot} (right) shows the second order convergence of $E_h$ and $E_h^\mathrm{pp}$ towards the reference energy and therefore also the asymptotical exactness of the lower bound. Note that, in general, the discrete energy alone is not a lower bound for the ground state energy, as the numerical example below for the constant potential shows.

\subsubsection*{Disorder potential and exponential localization}

\begin{figure}
	\includegraphics[height=.3\linewidth]{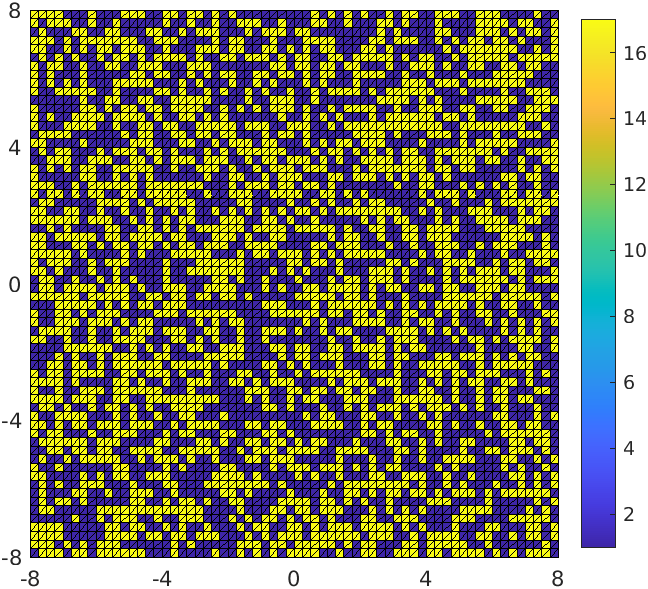}\hspace{1cm}
	\includegraphics[height=.3\linewidth]{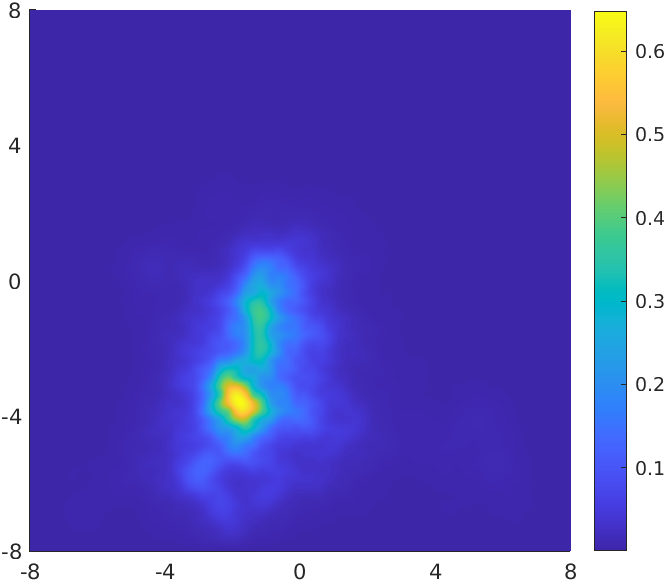}
	\caption{Disorder potential with coarse mesh used to construct the hierarchy of meshes in the background (left). Approximation of the highly localized ground state (right).}
	\label{fig:randpot}
\end{figure}

Second, we consider a disorder potential constructed using the Friedrichs-Keller triangulation shown in \cref{fig:randpot} (left). More precisely, we first join any pair of opposing triangles into squares of side length $\epsilon = 2^{-6}$. On all these squares, the coefficients is chosen to be constant, with values obtained as realizations of independent coin-flip random variables taking the values $1$ and $1+(2\epsilon L)^{-2}$. The parameter $\kappa$ is chosen to be one. For such coefficients there occurs an effect called Anderson localization (see, e.g., \cite{AltPV18,AltHP20,AltHP22} for numerical and theoretical studies). The exponential localization of the ground state can be seen in \cref{fig:randpot} (right). We emphasize that this example is numerically quite challenging, as can be seen from the comparatively large number of $J$-method iterations required.
For the discretisation, we use a hierarchy of meshes constructed by uniform refinement of the Friedrichs-Keller triangulation considered above. On each mesh of the hierarchy, the potential is obtained by prolongation.

\begin{figure}
	\includegraphics[width=.45\linewidth]{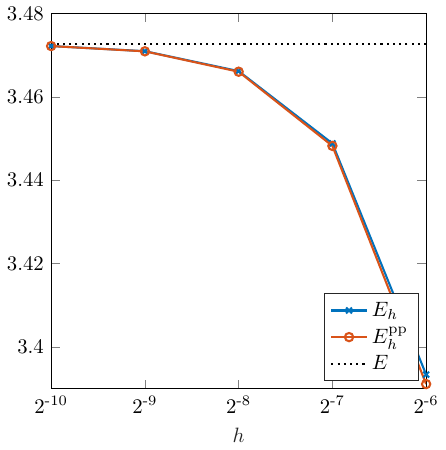}
	\hspace{.75cm}
	\includegraphics[width=.455\linewidth]{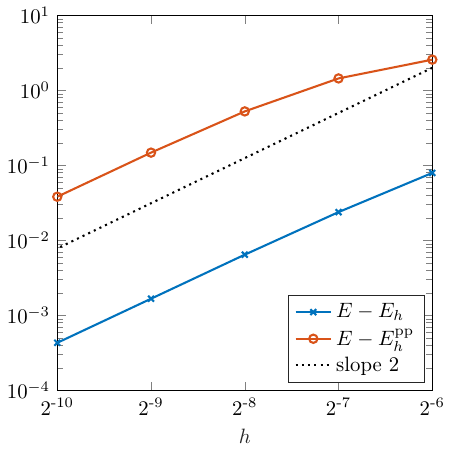}
	\caption{Ground state energy approximations for the disorder potential (left) and the difference between the reference energy and the energy approximations in a double-logarithmic plot (right). The blue and red curves correspond to the discrete energy and the post-processed discrete energy, respectively.}
	\label{fig:lowboundrandpot}
\end{figure}

In \cref{fig:lowboundrandpot} (left) it can be observed that, also for the disorder potential, $E_h$ and $E_h^\mathrm{pp}$ approach the ground state energy from below as the mesh size is decreased. \cref{fig:lowboundrandpot} (right) again demonstrates the second-order convergence of $E_h$ and $E_h^\mathrm{pp}$ towards the reference energy.

\subsubsection*{Constant potential and necessity of post-processing}
Third, we consider a constant potential, i.e., $V \equiv 1$. Although this choice may be unphysical, it is an example showing that the post-processing of the discrete energies is indeed necessary to obtain lower bounds. The parameter $\kappa$ is chosen to be one. For the discretisation we consider a hierarchy of meshes constructed by uniform refinement of the coarsest possible Friedrichs-Keller triangulation consisting of two elements.

\begin{figure}
	\includegraphics[width=.45\linewidth]{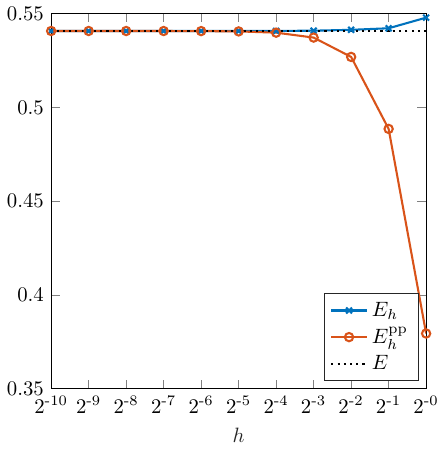}\hspace{.75cm}
	\includegraphics[width=.455\linewidth]{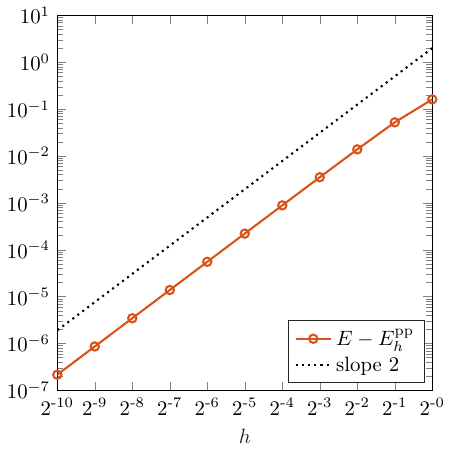}
	\caption{Ground state energy approximations for a constant potential (left) and the difference between the reference energy and the energy approximations in a double-logarithmic plot (right). The blue and red curves correspond to the discrete energy and the post-processed discrete energy, respectively.}
	\label{fig:lowboundconstpot}
\end{figure}

In \cref{fig:lowboundconstpot} (left), one observes that the discrete energies $E_h$ approach the ground state energy from above (and not from below) as the mesh size is decreased. Nevertheless, as predicted by \cref{th:lowerbound}, the post-processed energy $E_h^\mathrm{pp}$ is a lower bound. \cref{fig:lowboundconstpot} (right) shows the second-order convergence for $E_h^\mathrm{pp}$, while $E-E_h$ is negative in this example and therefore not shown in the double-logarithmic plot. 

\section{Conclusion}
In conclusion, this paper has effectively demonstrated the application of a mixed finite element discretisation to the Gross-Pitaevskii eigenvalue problem, with an emphasis on the computation of a lower energy bound. Our numerical experiments have not only validated the theoretical framework, but also confirmed the practicality of obtaining a computable lower bound on the ground state energy. This result provides a new aspect to the understanding and reliable numerical simulation of Bose-Einstein condensates.
\appendix
\section{Collection of frequently used bounds}
The following lemma provides, for any discrete function, 
a conforming lifting with a corresponding approximation estimate.

\begin{lemma}[Conforming lifting]\label{lem:appro:conti}
	For any $v_h\in \Uh$, there exists $v_h^c\in H^2(\Omega)\cap H_0^1(\Omega)$ such that it holds
	\begin{equation*}
		\tnof{G_hv_h-\nabla v_h^c} + \tnof{v_h-v_h^c}\lesssim h  \tnof{\ddiv G_hv_h}.
	\end{equation*}
\end{lemma}
\begin{proof}
	We denote by $v_h^c\in H^1_0(\Omega)$ the solution to Poisson's
	equation $-\Delta v_h^c = -\ddiv G_hv_h$ in $\Omega$ subject 
	to homogeneous Dirichlet boundary conditions. We emphasize that $v_h^c \in H^2(\Omega)\cap H^1_0(\Omega)$ with $\|v_h^c\|_{H^2}\lesssim \tnof{\ddiv G_h v_h}$, which follows from classical elliptic regularity theory on convex domains, see, e.g., \cite[Thm.~9.1.22]{Hac03}.  
	The pair $(v_h, G_h v_h)\in U_h\times \Sigma_h$ is the 
	Galerkin approximation of the mixed system and therefore satisfies the standard a priori error estimate 
	\begin{equation*}
		\tnof{G_hv_h-\nabla v_h^c} + \tnof{v_h-v_h^c}
		\lesssim h \|v_h^c\|_{H^2},
	\end{equation*}
cf.~\cite[Prop.~7.1.2]{BoffiBrezziFortin2013}. The assertion follows immediately.
\end{proof}

\begin{lemma}[$L^\infty$-bound]\label{lem:uh_Linfty}
Any $v_h\in \Uh$ satisfies that
\begin{equation*}
\|v_h\|_{L^\infty}\lesssim \|\ddiv G_h v_h\|_{L^2}.
\end{equation*}
\end{lemma}
\begin{proof}
 We denote by $v_h^c$ the conforming lifting and compute
 $$
 \|v_h\|_{L^\infty} 
  \lesssim
  \| v_h -\pi_h v_h^c\|_{L^\infty} + \|\pi_h v_h^c\|_{L^\infty}.
 $$
The first term on the right-hand side can be controlled by an
inverse estimate and the well-known superconvergence result
from \cite{DouglasRoberts1985,Brandts1994}. One obtains that
$$
\| v_h -\pi_h v_h^c\|_{L^\infty}
\lesssim h^{-d/2} \| v_h -\pi_h v_h^c\|_{L^2}
\lesssim h^{2-d/2} \|\ddiv G_h v_h \|_{L^2}.
$$
The remaining term is bounded by the $H^2$-norm of 
$v_h^c$, which again is controlled by $\|\ddiv G_h v_h \|_{L^2}$
thanks to elliptic regularity and the Sobolev embedding.
\end{proof}

Given $v_h\in U_h$, we define a piecewise affine function
$J v_h\in H^1_0(\Omega)$ by assigning to each vertex $z$ of the 
triangulation the arithmetic mean of the values that $v_h$
attains at $z$ when restricted to any elements containing $z$;
if $z$ is a boundary vertex, the value of $J v_h$ is set
to zero to conform to the homogeneous boundary condition.
Such averaging operators are well studied (see, e.g.,~\cite{BrS08})
and were used in the context of mixed finite elements, e.g.,
in \cite{HuangXu2012}.

\begin{lemma}[Averaging operator]\label{lem:properties_J}
Any $v_h\in \Uh$ satisfies that
\begin{equation*}
\|J v_h\|_{L^\infty}\lesssim \|v_h\|_{L^\infty}
\end{equation*}
and
\begin{equation*}
 \| h^{-1} ( v_h - J v_h)\|_{L^2} 
+\| \nabla J v_h\|_{L^2}
\lesssim  \| G_h v_h\|_{L^2}.
\end{equation*}
\end{lemma}
\begin{proof}
 The first bound follows directly from the construction of
 the function $J v_h$.
Following standard arguments, cf.~\cite[Lemma 10.6.6]{BrS08},
 we further obtain that
 $$
  \| h^{-1} ( v_h - J v_h)\|_{L^2} 
  +\| \nabla J v_h\|_{L^2}
  \lesssim
    \sqrt{
  \sum_{F}  h_F^{-1}\| [v_h]_F\|_{L^2(F)}^2
    },
 $$
where the sum runs over all faces $F$ and the bracket indicates
the inter-element jump across $F$, which is defined as the 
usual trace if $F$ is a boundary face.
 It was shown in \cite{LovadinaStenberg2006,GaoQiu2018} that this term 
 is bounded by $\| G_h v_h\|_{L^2}$.
\end{proof}

\begin{lemma}[Discrete embedding]\label{l:dem}
 Any $v_h\in U_h$ satisfies $\|v_h\|_{L^6}\lesssim \|G_h v_h\|_{L^2}$.
\end{lemma}
\begin{proof}
  Let $J v_h\in H^1_0(\Omega)$
 denote the regularization by averaging from above.
 From the triangle inequality, a classical comparison result between $L^p$-norms, and the Sobolev embedding, 
 we deduce that
 $$
  \|v_h\|_{L^6}
  \lesssim 
   \|h^{-d/3} (v_h - J v_h)\|_{L^2} + \|\nabla J v_h\|_{L^2}.
 $$
 By \cref{lem:properties_J} this is controlled by
 $\| G_h v_h\|_{L^2}$.
\end{proof}

As a consequence we note the following bound
\begin{equation}\label{eq:L6bound}
 \|u_h\|_{L^6} + \|u_h^c\|_{L^6}
 \lesssim \|G_h u_h\|_{L^2}+h \|\ddiv G_h u_h\|.
\end{equation}
\begin{proof}
 The bound for the first term on the left-hand side is
 shown in Lemma~\ref{l:dem}.
 Thanks to the Sobolev embedding, we have for the second term
 on the left-hand side that 
 $\|u_h^c\|_{L^6}\lesssim \|\nabla u_h^c\|_{L^2}$.
 With the triangle inequality and \cref{lem:appro:conti} we thus obtain that
 \begin{equation*}
 	\|u_h^c\|_{L^6} \lesssim \|G_h u_h\|_{L^2} + h \|\ddiv G_h u_h\|.\qedhere
 \end{equation*}
\end{proof}

\bibliographystyle{amsalpha}
\bibliography{bib}
\end{document}